\input amstex
\documentstyle{amsppt}
\voffset-,6in
\magnification=\magstep1
\NoBlackBoxes
\nologo
\topmatter
\title On certain one-counter shifts
\endtitle
\author  Wolfgang Krieger
\endauthor
\affil 
Institute for Applied Mathematics\\
University of Heidelberg\\
Im Neuenheimer Feld 294
69120 Heidelberg,
Germany\\
\endaffil
\abstract
Extrapolating from the two-block system of an example of a nonsofic shift hat was given by Lind and Marcus, a class of one-counter shifts is described, that is disjoint from the class of standard one-counter shifts.
\endabstract
\endtopmatter



\bigskip


\bigskip
\heading 1. Introduction
\endheading
Let $\Sigma$ be a finite alphabet. We use the notation 
$$
x_{[i,k]} = (x_i)_{i \leq j \leq k}, \qquad x \in \Sigma^{\Bbb Z},\quad i, k \in  \Bbb Z, i \leq k,
$$
also setting for $a= x_{[i.k]} $ and $ i \leq i^{\prime} \leq k^{\prime} \leq k$,
$$
a _{[i^{\prime},k^{\prime}]}   =   x _{[i^{\prime},k^{\prime}]}  .
$$
A block $x_{[i,k]}$ will also denote the word that is carries. On the shift space $ \Sigma^{\Bbb Z},$ there acts the shift by
$$
x \to (x_{i+1})_{i\in \Bbb Z},\quad x =  (x_{i})_{i\in \Bbb Z}, \in \Sigma^{\Bbb Z}.
$$
A closed shift invariant subset of $\Sigma^{\Bbb Z}$ is called a subshift. An introduction to the theory of subshifts is  in the book by Lind and Marcus \cite {LM}. A word is called admissible for a subshift if it appears in a point of the subshift. We denote the language of admissible words of a subshift $ X \subset \Sigma^{\Bbb Z}$   by $\Cal L(X)$, the length of a word $a \in  \Cal L(X) $ we denote by $\ell(a)$, and the set of $a \in  \Cal L(X) $ of length $n\in \Bbb N$ we denote by $ \Cal L_n(X) $. 
The  $n$-block system of a subshift $X \subset \Sigma^{\Bbb Z}$ is a subshift with alphabet  $\Cal L_n(X) $ 
and a topological conjugacy of $X$ onto its $n$-block system is  given by
$$
x \to (x_{[i, i+n)})^{}_{i\in \Bbb Z}  \ \ (x \in X).
$$

By using a matrix $(A(\sigma, \sigma^\prime))_{\sigma, \sigma^\prime \in \Sigma}$,
$$
A(\sigma, \sigma^\prime) \in\{  0, 1\} \quad \sigma, \sigma^\prime \in \Sigma,
$$
as a transition matrix one obtains a subshift $X_A$ by
$$
X_A = \{  (\sigma_{i}) _{i \in \Bbb Z} : A(\sigma_i, \sigma_{i+1} ) = 1, i \in \Bbb Z \}.
$$
These subshifts $X_A$ are called topological Markov shifts and they are special cases of subshifts of finite type $X_{\Cal F}$, where $\Cal F$ is a finite set of words in the alphabet $\Sigma$ (of length at least two), and where the subshift $X_{\Cal F}$ is  obtained by excluding the words in $\Cal F$  from  $\Sigma^{\Bbb Z}$. With $n $ the maximal length of a word in $\Cal F$ the $(n-1)$-block system of the subshift $X_{\Cal F}$ is a topological Markov shift. The subshifts of finite type belong to the class of sofic systems that are obtained from finite directed graphs, in which every vertex has at least one outgoing and at least one incoming edge, and that are labeled with symbols from the alphabet $\Sigma$, with the language of admissible words of the sofic system equal to the set of label sequences of finite paths in the graph.

The coded system \cite{BH} of a formal language ${\Cal C}$
in a finite alphabet
$\Sigma$
is the subshift
that is obtained as the closure of the set of points
in $\Sigma^{\Bbb Z}$
that carry bi-infinite concatenations of words in
${\Cal C}$.
${\Cal C}$ can here always be chosen 
to be a prefix code.
More generally, a Morkov code (see \cite{Ke}) 
is given by a formal language  ${\Cal C}$
of words in a finite alphabet $\Sigma$ together with
a  finite index set $\Gamma$, mappings 
$s: {\Cal C} \longrightarrow \Gamma,
 t: {\Cal C} \longrightarrow \Gamma
$
and a transition matrix
$(A(\gamma, \gamma'))_{\gamma,\gamma' \in \Gamma}, A(\gamma, \gamma')\in \{0,1\}, 
\gamma,\gamma' \in \Gamma.
$
From a Markov code $({\Cal C},s,t)$
one obtains its Markov coded system
 as the subshift that is the closure
of the set of points $x \in \Sigma^{\Bbb Z}$
such that there are indices 
$i_k \in {\Bbb Z}, k \in {\Bbb Z}$,
$i_k < i_{k+1}, k \in {\Bbb Z},$
such that 
$x_{[i_k,i_{k+1})} \in {\Cal C}, k \in {\Bbb Z}$,
and such that 
$$
A(t(x_{[i_{k-1},i_{k})}), s(x_{[i_k,i_{k+1})}) ) = 1, \qquad k \in {\Bbb Z}.
$$

Given a subshift $X \subset {\Sigma}^{\Bbb X}$
a word $v \in {\Cal L}(X)$
is called synchronizing if for 
$u, w \in {\Cal L}(X)$ 
such that 
$ uv, vw \in {\Cal L}(X)$
also
$ u v w \in {\Cal L}(X)$.
A topologically  transitive
subshift is called synchronizing 
if it has a synchronizing word. As in \cite {KM} we say that a synchronizing subshift $X \subset \Sigma^{\Bbb Z}$ is strongly synchronizing
if there exists a $Q \in \Bbb Z_+$,  such that the following holds:
if $ x \in X$ and $I_{-}, I_{+} \in {\Bbb Z}, I_{-} < I_{+}$
are such that  
$x_{[I_{-}, I_{+}]}$ is synchronizing, then 
there exists an index $i,$  $I_{-} - Q \le i \le I_{+} + Q,$
such that 
$x_{i}$ is a synchronizing symbol.
Denote for a strongly synchronizing subshift $X \subset \Sigma^{\Bbb Z}$  its set of synchronizing 
symbols by $\Sigma_{synchro}(X)$ 
 and by $\Cal B (X)$ the set of words in $\Cal L(X)$ of length at 
least two that begin and end with a synchronizing symbol with no synchronizing symbol in between.
$\Cal B (X)$ determines Markov codes $\Cal C ^-(X)$ and $\Cal C ^+(X)$ that have as their index set $\Gamma$ a set of subsets of  $\Sigma_{synchro}(X)$. The Markov code $\Cal C ^-(X)$ contains the words that are obtained by removing from the words in $\Cal B (X)$ their last symbols, and  for a word $c\in \Cal C ^-(X)$ the
set $t(c)$ is equal to the set of synchronizing symbols
that can follow $c$ and the
set $s(c)$ is equal to the singleton set  that contains the first symbol of $c$, and the positive entries of the  transition matrix of the  Markov code are given by
$$
A({\Sigma_\circ}, \{ \sigma\}) = 1, \qquad
\Sigma_\circ \in \bigl\{ \{ t(c) \mid c \in {\Cal C}^-(X)\bigr \}\},\quad
 \sigma \in \Sigma_\circ, \quad
\sigma \in \Sigma_{synchro}(X) .
$$
$\Cal C ^+(X)$ is symmetric to $\Cal C ^-(X)$. A strongly synchronizing subshift $X$ is the Markov coded system of $\Cal C ^-(X)$, as well as of $\Cal C ^+(X)$, and it is seen  that $X$ can be reconstructed from $\Cal B (X)$.

A formal language is said to be a one-counter language if its words  can be recognized by a push-down automaton with one stack symbol.
We say that a strongly synchronizing subshift $X$ is a one-counter shift,  if   $\Cal B (X)$ (or $\Cal C ^-(X)$, or  $\Cal C ^+(X) $) as a formal language is a one-counter language.
In \cite {KM} the  class of  standard one-counter shifts $X \subset \Sigma^{\Bbb Z}$ was introduced by specifying structural properties of  $\Cal C ^-(X)$.  With the alphabet $\{a, b, c\}$ prototypical examples of standard one-counter shifts are the coded systems of the code  
$$
\Cal C_{reset} = \{ ab^kc^l: k, l \in  \Bbb N, l \leq k \},
$$
and of the code
$$
\Cal C_{counter} = \{ ab^kc^l: k, l \in  \Bbb N, l  =  k \}.
$$
The coded system of $\Cal C_{reset} $ is a prototypical example of what was called in \cite{KM}
 a one-counter shift with reset, and the coded system of $\Cal C_{counter} $ is a prototypical example of what was called in \cite{KM}
 a one-counter shift without reset.
The class of standard one-conter shifts with reset is not closed under taking inverses. 

In the present paper our starting point is the two-block system of a subshift that appeared as an example of a nonsofic subshift in the book by Lind and Marcus (\cite {LM}, Example 1.2.9). We introduce this subshift, that we call the Lind-Marcus shift, as the subshift with alphabet $\{a, b, c\}$ that does not allow the words
$
ab^k c^la, \thinspace k, l \in \Bbb Z_+, k  \ne  l
$
(as in \cite {Ma1}, see also \cite{Ma2} section 3). In other words, the Lind-Marcus shift is the coded system of the code that contains besides the words $ab^kc^l, k \in \Bbb Z_+ ,$  the words that begin with a symbol $a$ that is then followed by a word in the symbols $b, c$ that is not equal to any of the words $b^kc^l, k, l \in \Bbb Z_+, k \ne l.$
In  \cite {M1} Kengo Matsumoto has computed the K-groups of the Lind-Marcus shift. He also showed that the C*-algebra that is associated to the Lind-Marcus shift is simple and purely infinite. 

In section 2 we look at the two-block system of the Lind-Marcus shift. In section 3, abstracting  properties of the two-block system of the Lind-Marcus shift, we describe a class of one-counter shifts  that is disjoint from the class of standard one-counter shifts.
We call the subshifts in this class  Lind-Marcus one-counter shifts. The class of Lind-Marcus one-counter shifts  is closed under taking inverses. Because of the time symmetry we find it appropriate to identify the class of Lind-Marcus one-counter shifts $X \subset \Sigma^{\Bbb Z}$ by specifying structural properties of $\Cal B(X)$, rather than of   $\Cal C ^-(X)$ (or of  $\Cal C ^+(X)$).

Recall that, given subshifts $X \subset \Sigma^{\Bbb Z}$,
 $\tilde{X} \subset \tilde{\Sigma}^{\Bbb Z}$,
 and a topological conjugacy 
 $\tilde\varphi:\tilde X \to X$, 
 there is for some $L \in \Bbb Z_+$ a block map
 $\tilde\varPhi:\tilde X_{[-L,L]} \to \Sigma$,
 such that 
 $$
\tilde\varphi(\tilde x) = (\tilde\Phi(\tilde x_{[i-L,i+L]}))_{i \in {\Bbb Z}}, \qquad \tilde x \in \tilde X.
$$
We will use the notation
$$
\tilde\Phi( \tilde a) = ( \tilde\Phi (  \tilde a  _{[i-L, i + L]} )_{I_--L\leq i \leq I_++L }, \quad  \tilde a \in \tilde X_{[I_-,I_+]},I_-,I_+ \in \Bbb Z, I_-\leq I_+.
$$

In section 4 we show that in a situation where one is given subshifts 
$
X \subset \Sigma^{\Bbb Z},
\tilde{X} \subset \tilde{\Sigma}^{\Bbb Z},
$
 $\tilde X$ a Lind-Marcus one-counter shift, 
 and a topological conjugacy $\tilde\varphi:\tilde X \to{X} $ that is given by a one-block map
$\tilde\Phi:\tilde\Sigma \to{\Sigma}$,
such that 
$$
\tilde\Phi^{-1}({\Sigma}_{synchro}({X}) \subset \Sigma_{synchro}(\tilde X),
$$
one has that $\tilde{X} $ is also a  Lind-Marcus one-counter shift. This leads one to subshifts of Lind-Marcus one-counter type, that one introduces as the subshifts that have an $n$-block system that is a Lind-Marcus  one-counter shift, and to the result that a subshift that is topologically conjugate to a subshift of Lind-Marcus one-counter type is itself a subshift of Lind-Marcus one-counter type. Compare here the relationship between topological Markov shifts and subshifts of finite type (also see  Section 3c of \cite{KM}). The proof proof of  simplicity and pure infinitedness of the assciated C*-algebra \cite{M1, M2} carries over to the  subshifts of Lind-Marcus one-counter type.

\heading 2. The two-block system of the Lind-Marcus shift 
\endheading

Let $Y$ be the two-block system of the Lind-Marcus shift. A word $(u_i)_{1 \leq i \leq I}, I > 1,$ that is admissible for the Lind-Marcus shift, determines a word $(u)^{\langle2\rangle}$ that is admisssible for $Y$ by
$$
(u)^{\langle2\rangle} = (u_iu_{i+1})_{1 \leq i < I}.
$$
\proclaim{Lemma 2.1}
$$
\align
\Cal B (Y) =
&
 \{ (ab^kc^ka)^{\langle2\rangle}: k \in \Bbb Z_+\} \cup \{ (cb^kc^lb)^{\langle2\rangle}: k, l \in \Bbb Z_+\}\cup
  \\
 &
 \{ (cb^kc^la)^{\langle2\rangle}: k, l \in \Bbb N\}\cup  \{ (cb^kc^lb)^{\langle2\rangle}: k \in \Bbb N\} \cup
 \\
 &
\{ (ac^lb)^{\langle2\rangle}: l \in \Bbb Z_+\} \cup
\{ (cb^ka)^{\langle2\rangle}: k \in \Bbb Z_+\}  .
 \endalign
$$
\endproclaim
\demo{Proof}
The synchronizing symbols of $Y$ are
$$
aa,\quad ab, \quad ac, \quad ba,\quad ca, \quad cb,
$$
and the set of words in $\Cal L(Y)$ that do not contain a ysnchronizing symbol is equal to
$$
\{(b^kc^l)^{\langle2\rangle}  : k, l \in \Bbb Z_+, k+l > 1  \}.
\qed
$$
\enddemo
\proclaim{Lemma 2.2} 
$Y$ is strongly synchronizing.
\endproclaim
\demo{Proof}
The set of non-synchronizing words in $\Cal L(Y)$ coincides with the set of words in $\Cal L(Y)$ that do not contain a synchronizing symbol.
\qed
\enddemo
\heading 3. Lind-Marcus one-counter shifts
\endheading

Given a subshift $X \subset \Sigma^{\Bbb Z}$ we set
$$
X_{[i,k]} = \{ x_{[i,k]}  : x \in X \}, \quad i, k \in \Bbb Z, i \leq k.
$$
We also set
$$
\Gamma^-_n (a)= \{b \in \Cal L_n (X): ba \in   \Cal L (X)   \}, \quad n \in \Bbb N,
$$
$$
\Gamma^-(a) = \bigcup_{n \in \Bbb N}\Gamma^-_n (a), \quad a \in  \Cal L (X) .
$$
$\Gamma^+$ has the symmetric meaning.

We recall some notions and results from \cite{KM}.
Let $X$ be a topologically transitive subshift.
We say that a pair  
$((\alpha_{-})_{i\in {\Bbb Z}},
(\alpha_{+})_{i\in {\Bbb Z}})
$
of fixed points of $X$ is a characteristic pair of fixed points of $X$
if it is the unique pair of fixed points of $X$ that satisfies the following conditions $(a), (b) $ and $(c^-)$
and a condition $(c^+)$
that is symmetric to $(c^-)$:

$(a)$ There is an orbit $O_X$ in $X$ that contains all points that are 
left asymptotic to $(\alpha_{-})_{i\in {\Bbb Z}}$
and 
right asymptotic to $(\alpha_{+})_{i\in {\Bbb Z}}$
and that do not contain a synchronizing word.

$(b)$ $X$ has a point that is
left asymptotic to $(\alpha_{+})_{i\in {\Bbb Z}}$
and 
right asymptotic to $(\alpha_{-})_{i\in {\Bbb Z}}$
and that contains a synchronizing word.

$(c^-)$  There exists a $K \in {\Bbb N}$ 
such that the following holds:
If $ x \in X$ and $I_{-}, I_{+} \in {\Bbb Z}, I_{-} \le I_{+}$
are such that  
$x$ is
right asymptotic to $(\alpha_-)_{i\in {\Bbb Z}}$
and 
$x_{[I_{-}, I_{+}]}$ is synchronizing, 
and
$x_{(I_+, I_+ +k]}, k \in {\Bbb N}$
is not synchronizing,
then 
there exists an index $i,$  $I_{-} < i \le I_{+} + K$
such that 
$x_{j} = \alpha_{-}, j \ge i$. 

If $((\alpha_{-})_{i\in {\Bbb Z}},
(\alpha_{+})_{i\in {\Bbb Z}})$ is a characteristic pair of fixed point of $X$, and if $\varphi$ is a topological conjugacy of a subshift $\tilde X$ onto $X$ then 
 $
 (\varphi^{-1}
 ((\alpha_{-})_{i\in \Bbb Z},\varphi^{-1}
 ((\alpha_{+})_{i\in \Bbb Z}) $
is a characteristic fixed point of $\tilde X$.

Given a fixed point $(\alpha)_{i \in {\Bbb Z}}$
of a subshift $X \subset {\Sigma}^{\Bbb Z}$,
and given a symbol $\sigma \in \Sigma_{synchro}$
we denote by ${\Cal D}(\sigma,\alpha)$ 
the set of words $d^-$ such that 
$$
\sigma d^- \in \Gamma^-(\alpha^k), \qquad k \in {\Bbb N},
$$
and such that in case that 
$d^-$ is not empty,
$d^-$ does not contain a synchronizing symbol and ends in a symbol that is different from $\alpha$.
A set 
${\Cal D}(\alpha,\sigma)$
of words is defined symmetrically.
We set for a strongly synchronizing subshift $X \subset \Sigma^{\Bbb Z}$
with a characteristic pair 
$(\alpha_{-})_{i\in {\Bbb Z}}, 
(\alpha_{+})_{i\in {\Bbb Z}}
$
of fixed points
$$
\align
\Sigma_{-}(X) & = 
\{ \sigma_{-} \in \Sigma_{synchro} : {\Cal D}(\sigma_{-},\alpha_{-}) \ne \emptyset \},\\
\Sigma_{-}^+(X) & = 
\{ \sigma_{-} \in \Sigma_{synchro}: {\Cal D}(\sigma_{-},\alpha_{+}) \ne \emptyset \},\\
\Sigma_+(X) & = 
\{ \sigma_{+} \in \Sigma_{synchro} : {\Cal D}(\alpha_{+},\sigma_+) \ne \emptyset \},\\
\Sigma_{+}^-(X) & = 
\{ \sigma_{+} \in  \Sigma_{synchro}:{ \Cal D}(\alpha_{-},\sigma_{+}) \ne \emptyset \}.
 \endalign
$$
The sets $\Sigma_{-}(X)$ and $\Sigma_{+}(X)$
are not empty and the sets 
$\Cal D(\sigma_{-},\alpha_{-}), \, \sigma_{-}\in \Sigma_{-}(X)$
and
${\Cal D}(\alpha_{+},\sigma_{+}), \, \sigma_{+}\in \Sigma_{+}(X)$
are finite.

Let $X \subset \Sigma^{\Bbb Z}$ be a subshift 
with a characteristic pair 
$(\alpha_{-})_{i\in {\Bbb Z}}, 
(\alpha_{+})_{i\in {\Bbb Z}}
$
of fixed points.
Let $x \in O_X$.
If for some $i_{\circ} \in {\Bbb Z}$,
$$
\align
x_i & = \alpha_-, \qquad i \le i_{\circ},\\
x_i & = \alpha_+, \qquad i > i_{\circ}
\endalign
$$
then set $c_X$ equal to the empty word.
Otherwise determine 
$i_-, i_+ \in {\Bbb Z}, i_- < i_+$, by
$$
\align
x_i & = \alpha_-, \qquad i < i_-,\\
x_{i_-}& \ne \alpha_-, \\
x_{i_+}&  \ne \alpha_+,\\
x_i & = \alpha_+, \qquad i > i_+
\endalign
$$
and set 
$c_X$ equal to the word 
$x_{[i_-,i_+]}$.

In preparation of the introduction of the Lind-Marcus one-counter shifts we formulate for
a strongly synchronizing subshift with a characteristic pair of fixed points $((\alpha_{-})_{i\in {\Bbb Z}},  
$ $ (\alpha_{+})_{i\in {\Bbb Z}})
$ two conditions $(R^-)$ and $(R^-)$ that are symmetric to one another. 
 Conditon  $(R^-)$ is as follows:

$(R^-)$ There exists an $R^-\geq \ell (c_X)$, such that, if 
$$
\sigma_-^+ \in \Sigma_-^+(X)  , d^- _+\in \Cal D(\sigma_-^+, \alpha_+), \ell(d^-_+) > R^-  ,
$$
then $ \sigma_-^+ \in \Sigma_-(X)$ and there is a $d^- \in \Cal D(\sigma_-^+, \alpha_-)$ such that
$$
 d^- _+ = d_-\alpha_-^{\ell( d^- _+ ) - \ell ( d^- ) - \ell (c_X)}c_X.
$$
For a subshift $X$ that satisfies condition $(R^-)$  we denote by $R^-(X)$ the smallest $R^-$
 such that $(R^-)$ holds. $R^+(X)$ is defined symmetrically.


For a strongly synchronizing subshift  $X \subset \Sigma^{\Bbb Z}$  
with a characteristic pair 
$((\alpha_{-})_{i\in {\Bbb Z}}, $
$(\alpha_{+})_{i\in {\Bbb Z}})
$
of fixed points, that satisfies conditions $(R^-)$ and $(R^+)$, we denote by 
$
\Xi_-(X) 
$
$
(
\Xi_-^+(X) 
)
$
the set of words $ \sigma_- d^- , \sigma_- \in \Sigma_-(X), d^- \in \Cal D(\sigma_- , \alpha_{-}  )$
($ \sigma_-^+ d^-_+ , \sigma_-^+ \in \Sigma_-^+(X), d^-_+ \in \Cal D(\sigma_- ^+, \alpha_+  )$), such that there is a 
$K_- \in \Bbb N$
($
K_-^+\in \Bbb N
$)
such that  for $  \sigma_+^- \in  \Sigma_+^-(X), d^+_- \in \Cal D(\alpha_{-} ,  \sigma_+^-   )$
$
(  \sigma_+ \in  \Sigma_+(X), d^+ \in \Cal D( \alpha_+,  \sigma_+   ))
$, and for 
 $k_- \geq  K_-$ 
 $
 ( k_-^+ \geq  K_-^+    )
 $
  the word
 $ \sigma_- d^-\alpha_-^{k_- }d^-_+  \sigma_-^-   $
 $
 (    \sigma_-^+ d^-_+  \alpha_+^{k_-^+}  d^+   \sigma_+  )
 $
  is admissible for $X$. Also, in the case that $ \Xi_-(X) \neq \emptyset  $  we denote by $K^-(X) $($K^+_-(X) $) the smallest $ K^- \in \Bbb N$ ($K_-^+  \in \Bbb N$) such that 
  $$
   \sigma_- d^-\alpha_-^{k_- }d^-_+  \sigma_-^-  \in \Cal L (X), \qquad k_-  > K^-,  \sigma_- d^- \in  \Xi_-(X) 
  $$
  $$
  (  \sigma_-^+ d^-_+  \alpha_+^{k_+}  d^+   \sigma_+  \in \Cal L (X),\qquad k_- ^+ > K_- ^+ ,  \sigma_-^+ d^-_+ \in  \Xi_-^+(X) 
  ),
  $$
  where we observe, that by condition $(R^+)$,
  $$
  \{   \sigma_- d^-  \alpha_-^{k_- } c_X  : k_-  > K_- \}\subset \Xi_-^+(X) . \tag 1
  $$
$\Xi_+(X) ,\Xi_+^-(X) K_+(X),K_+^-(X)$ have the symmetric meaning. We will find it necessary that for long words the converse inclusions to (1) and to its symmetric counterpart hold, and in order to ensure this, we impose
on a strongly synchronizing subshift  $X \subset \Sigma^{\Bbb Z}$  
with a characteristic pair 
$((\alpha_{-})_{i\in {\Bbb Z}}, $
$(\alpha_{+})_{i\in {\Bbb Z}})
$
of fixed points, that satisfies conditions $(R^-)$ and $(R^+)$,
two  conditions $(R_\Xi^-)$ and  $(R_\Xi^+)$, that are symmetric to one another. $(R_\Xi^-)$ is as follows:

 $(R_\Xi^-)$
There exists an $R^-_\Xi \geq R^-(X)$,
 such that, if $$  \sigma_-^+ d^-_+ \in  \Xi_-^+(X)    , \ell (  d^-_+  ) >  R^-_\Xi  ,  $$ then $ \sigma_-^+   \in  \Sigma_-(X)$ and there is a $d^- \in \Cal D(  \sigma_-^+   , \alpha_-) $ such that  $\sigma_-^+ d^-  \in  \Xi_-(X) $, and 
 $$
d^-_+ = d^-  \alpha_-^{\ell ( d^-_+  )- \ell ( d^- ) -\ell ( c_X ) }  c_X .
 $$
 For a subshift $X$ that satisfies conditions $(R_\Xi^-)$ we denote by $R^-_\Xi(X)$ the smallest
 such that $(R_\Xi^-)$ holds. $R^-_\Xi(X)$ is defined symmetrically.

For a strongly synchronizing subshift  $X \subset \Sigma^{\Bbb Z}$  
that has a characteristic pair 
of fixed points, and that satisfies conditions $(R^-),(R^+),(R_\Xi^-),(R_\Xi^+)$, we set
$$
\multline
\Cal B^{(X)} (\Xi, \alpha_ - , \alpha_ +)=\\
\{ \sigma_-d^- \alpha_-^{k_-}c_X\alpha_+^{k_+}d^+\sigma_+:
\sigma_-d^- \in \Xi_-(X), 
\sigma_+ \in \Sigma_+(X) ,d^+\in \Cal D(\alpha_+ ,  \sigma_+ ),\\
k_-,k_+ \in \Bbb N,  \ell (d^-  ) + k_ -+ \ell( c_X ) > R^-_\Xi(X) , \ell( c_X )+  k_ -+  \ell (d^- ) > R^+_\Xi(X)\}.
\endmultline
$$
$\Cal B^{(X)} ( \alpha_ - , \alpha_ +, \Xi)$ is defined symmetrically. Also set
$$
\align
\Cal B^{(X)} (\Xi, \alpha_ -)=
&\{ \sigma_-d^- \alpha_-^{k_-}d_-^+\sigma_+^-: \\ 
&\sigma_-d^- \in \Xi_-(X),
\sigma_+ \in \Sigma_+(X) ,d^+\in \Cal D(\alpha_+ ,  \sigma_+ ),\\
&k_- \in \Bbb N,  \ell (d^-  ) + k_ -+ \ell( c_X ) > R^-_\Xi(X), \ell (d_-^+) \leq R^+(X)\}.\\
\endalign
$$
and
$$
\align
\Cal B^{(X)} (\Xi, \alpha_ +)=
\{ &\sigma_-^+d^- _+\alpha_+^{k_+}d^+\sigma_+:\\
&\sigma_-^+d^- _+ \in \Xi_-^+(X), \sigma_+ \in \Sigma_+(X) ,d^+ \in \Cal D( \alpha_+,\sigma_+ ),\\
&\ell (d_-^+) 
\leq R^- _\Xi(X), k_+ \in \Bbb N,  \ell( c_X ) + k_ + + \ell (d^+ )  > R^+(X) \}.
\endalign
$$
$\Cal B^{(X)} (\alpha_ +, \Xi) $ and $\Cal B^{(X)} ( \alpha_ -, \Xi)$  are defined symmetrically.
Also, given $J_-, J_+ \in \Bbb Z_+$ and mappings 
$$
(\sigma_-,d^-) \to \Delta^-(\sigma_-, d^-)\subset \Bbb Z_, \quad 
\sigma_- \in \Sigma_-(X), d^-\in {\Cal D}^-(\sigma_-,\alpha_-),\sigma_-d^- \notin \Xi_-(X),
$$
$$
(d^+, \sigma_+) \to \Delta^+(d^+, \sigma_+)\subset  \Bbb Z_+ \quad 
\sigma_+ \in \Sigma_+(X), d^+\in {\Cal D}^-(\alpha_+, \sigma_+),
d^+ \sigma_+\not\in \Xi_+(X),
$$ 
set
$$
\align
 \Cal B_{counter}^{(X)}(&\Delta^-, J_-, J_+,\Delta^+)  =  \\
\{ &\sigma_-d^- \alpha_-^{k_-}c_X\alpha_+^{k_+}d^+\sigma_+:
\\
&
\sigma_- \in \Sigma_-(X) ,d^- \in \Cal D(\sigma_-  ,   \alpha_-), \sigma_-d^- \notin \Xi_-(X),
\\
&
\sigma_+ \in \Sigma_+(X) ,d^+\in \Cal D(\alpha_+ ,  \sigma_+ ), d^+\sigma_+ \notin \Xi_+(X),
\\
&
k_-,k_+ \in \Bbb N,  \ell (d^-  ) + k_ -+ \ell( c_X ) > R^-(X) , \ell( c_X )+  k_ ++  \ell (d^- ) > R^+(X), \\
&( \Delta^-(\sigma_-, d^-) + k_- + J_- ) \cap (J_+ + k_+ + \Delta^+(d^+,\sigma_+)) \ne \emptyset  \}.
\endalign
$$

We define a Lind-Marcus one-counter shift as a strongly synchronizing subshift with a characteristic pair of fixed points,  that satisfies conditions $(R^-)$,$(R^+)$,$(R^-_\Xi)$,
$(R^+_\Xi)$, and that is such that  
$$
\Xi_-(X)\neq \emptyset, \quad
\{   \sigma_-d^-:   \sigma_- \in \Sigma_-(X), d^- \in \Cal D (  \sigma_-  ,  \alpha_-   )  \} \setminus  \Xi_-(X)
\neq \emptyset,
$$
$$
\Xi_+(X) \ne \emptyset, \quad
\{  d^+ \sigma_+:   \sigma_+ \in \Sigma_+(X) ,d^+\in \Cal D (  \alpha_+ , \sigma_+    )  \}
 \setminus \Xi_+(X)
\neq \emptyset,
$$
 and such that there are
 $J_-, J_+ \in \Bbb Z_+$ and mappings 
$$
(\sigma_-,d^-) \to \Delta^-(\sigma_-, d^-)\subset \Bbb Z_, \quad 
\sigma_- \in \Sigma_-(X), d^-\in {\Cal D}^-(\sigma_-,\alpha_-),\sigma_-d^- \notin \Xi_-(X),
$$
$$
(d^+, \sigma_+) \to \Delta^+(d^+, \sigma_+)\subset  \Bbb Z_+ \quad 
\sigma_+ \in \Sigma_+(X), d^+\in {\Cal D}^-(\alpha_+, \sigma_+),
d^+ \sigma_+\not\in \Xi_+(X),
$$ 
and an $I \in \Bbb N$ such that
$$
\align
\{b \in \Cal B(X): \ell(b) > I\} = \  &   
\{b \in  {\Cal B}_{counter}^{(X)}(\Delta^-, J_-, J_+,\Delta^+) \  \cup  \tag LM\\
 &\Cal B^{(X)} (\Xi, \alpha_ - , \alpha_ +) \cup \Cal B^{(X)} ( \alpha_ - , \alpha_ +,\Xi) \ \cup \\
&\Cal B^{(X)} (\Xi, \alpha_ -)  \cup \Cal B^{(X)} ( \alpha_ - ,\Xi) \ \cup\\
&\Cal B^{(X)} (\Xi,  \alpha_ +) \cup \Cal B^{(X)} ( \alpha_ +,\Xi): \ell(b) > I\}.   
\endalign
$$
If (LM) holds then we say that $J_-, J_+ , \Delta^-, \Delta^+$ are parameters for the Lind-Marcus shift $X$.

$Y$ is a Lind-Marcus one-counter shift. One has
$$
\Sigma_{-} (Y) = \{ ab, cb  \}    ,  \  \  \Sigma_{-} ^{+}(Y) = \{ ac  \}  ,   
$$
$$
\Sigma_{+} (Y) =  \{ ca, cb   \}  ,   \  \  \Sigma_{+}^{-} (Y) = \{ba   \},
$$
the sets $\Cal D(ab  , bb ),  \Cal D( cb , bb ), \Cal D( ac , cc ) $ and the sets $\Cal D(cc  , ca ),  \Cal D( cc ,  cb) , \Cal D( bb , ba )$  are here equal to the sigleton set that contains the empty word, and
$$
\align
&\Xi_{-} (Y) = \Xi_{+} (Y) = \{ cb  \}    ,  \\  
&\Xi_{-} ^{+}(Y) = \{ (ab^kc)^{\langle 2 \rangle} : k \in \Bbb Z_+ \}  , \  \  \Xi_{+}^{-} (Y) = \{(bc^ka)^{\langle 2 \rangle} : k \in \Bbb Z_+  \}  .
\endalign
$$
From Lemma (2.1) it can be seen that (LM) holds for $Y$.

With the alphabet $\{   b, c\} \cup\{a_n: 1 \leq n \leq N\},N> 1, $ one has Lind-Marcus one-counter shifts, that are closely patterned after the Lind-Marcus shift, and  the 2-block systems of
that do not allow the words
$
a_nb^k c^la_m, \quad k, l \in \Bbb Z_+, k  \ne  l, 1 \leq n, m   \leq N.
$
In  \cite {M2} Kengo Matsumoto has computed the K-groups of 
these Lind-Marcus one-counter shifts.

For another example of a Lind-Marcus counter shift that is patterned after the 2-block system of the Lind-Marcus shift, take the alphabet $\{a, b, c, d\}$, is the synchronizing shift $X$ such that 
$$
\Cal C^-(X) = \{ ab^kc^k  : k \in \Bbb Z_+\} \cup       \{ ab^kc^l  : k,l \in \Bbb Z_+\} \cup 
\{ db^kc^l  : k, l \in \Bbb Z_+\} ,
$$
with
$$
\align
&t(  ab^kc^k ) = \{  a \} ,\quad k \in \Bbb Z_+, \\
&t(  ab^kc^l ) = \{  a, d \} ,\quad k \in \Bbb Z_+, k \neq l,  \\
& t(  db^kc^l ) = \{  a, d \} \quad k \in \Bbb Z_+.
\endalign
$$

\heading 4. Topological conjugacy
\endheading


For the proof of the next lemma compare the proof of  Lemma 3.10 in  \cite {KM}. We will  begin the proof by recalling  a construction from  Section 3 b of \cite {KM}.

\proclaim{Lemma 4.1} 
Let there be
given subshifts 
$
X \subset \Sigma^{\Bbb Z},
\tilde{X} \subset \tilde{\Sigma}^{\Bbb Z},
$ 
and a topological conjugacy $\tilde\varphi : \tilde X \to X $ that is given by a one-block map
$\tilde\varPhi:\tilde\Sigma \to {\Sigma}$,
such that 
$$
\tilde\varPhi^{-1}(\Sigma_{synchro}(X) \subset \Sigma_{synchro}(\tilde X).
$$
Let $X$ be a Lind-Marcus one-counter shift. Then $\tilde{X}$ is also a  Lind-Marcus one-counter shift.
\endproclaim
\demo{Proof}
Let  the inverse of $\tilde\varphi$ be given by for some $L \in \Bbb Z_+$ by a block map
 $\varPhi:X_{[-L,L]} \to \tilde{\Sigma}$. Let  $(({\alpha}_-)_{i \in {\Bbb Z}},
({\alpha}_+)_{i \in {\Bbb Z}})
$ be the characteristic pair of fixed points of $X$, and
set
$
(\tilde{\alpha}_-)_{i \in {\Bbb Z}} = \tilde{\varphi}^{-1}(
(\alpha_-)_{i \in {\Bbb Z}})
$,
$(\tilde{\alpha}_+)_{i \in {\Bbb Z}} = \tilde{\varphi}^{-1}(
(\alpha_+)_{i \in {\Bbb Z}}).$
Also choose a
$Q \in {\Bbb N}$ such that 
for a synchronizing word $a$ of $X$ and for 
$a^- \in \Gamma_Q^-(a), 
a^+ \in \Gamma_Q^+(a)$
the word 
$a^- a a^+$ contains a synchronizing symbol.

For
$\tilde{\sigma}_- \in \tilde{\Sigma}_-(\tilde{X}),
\tilde{b}^-\in \Gamma_{Q + L}^-(\tilde{\sigma}_-),$
one can set 
$$
\tilde{\varPhi}(\tilde{b}^-\tilde{\sigma}_-) = 
b^-(\tilde{b}^-\tilde{\sigma}_-) \sigma_-(\tilde{b}^-\tilde{\sigma}_-) a^-(\tilde{b}^-\tilde{\sigma}_-),
$$
where the words 
$b^-(\tilde{b}^-\tilde{\sigma}_-) $
and 
$a^-(\tilde{b}^-\tilde{\sigma}_-) $
and the symbol
$\sigma_-(\tilde{b}^-\tilde{\sigma}_-) $
are uniquely determined by
$\tilde{b}^-\tilde{\sigma}_-$
 under the condition that 
$\sigma_-(\tilde{b}^-\tilde{\sigma}_-) $
is synchronizing and that 
$a^-(\tilde{b}^-\tilde{\sigma}_-) $
does not contain a synchronizing symbol.
We set
$$
I_-(\tilde{b}^-\tilde{\sigma}_-) = \ell(a^-(\tilde{b}^-\tilde{\sigma}_-)),
$$
and for $\tilde d^- \in \Cal D( \tilde{\sigma}_- ,\tilde \alpha _-)$ we denote by $d^-(  \tilde{b}^-\tilde{\sigma}_-\tilde d ^-)  $ the longest prefix of the word  $a^-(\tilde{b}^-\tilde{\sigma}_-)\tilde\Phi ( \tilde d^-  )\alpha_-$ that is in $\Cal D(\sigma_-(\tilde{b}^-\tilde{\sigma}_-) , \alpha_-)$.

A converse construction yields for given  $\sigma_-  \in \Sigma _-(X)$ and
 $d^- \in \Cal D( \sigma_-  , \alpha_-)$
 a $\tilde \sigma_- \in \Sigma_-(\tilde X)$ and 
 $\tilde b_-\in \Gamma^-_{Q+L}(  \tilde \sigma_-   ), 
 \tilde d^-
 \in \Cal D(  \tilde \sigma_-   , \alpha_-)$ such that 
$$
\sigma_- = \sigma_-(\tilde b_-  \tilde \sigma_-       ), 
\quad d^- =   (\tilde b_-  \tilde \sigma_-    \tilde d^-      ) . \tag 1 
$$
With a $c^- \in \Gamma^-_{Q+2L}(  \sigma_-   ),$ one lets  $  \tilde \sigma_- $ be the last synchronizing symbol in the word 
$\Phi(c^-  \sigma_-  d^- \alpha_- ^{2L + 1}  )$ and one lets $\tilde d^- (  \tilde \sigma_-   ,   \tilde \alpha_-  )$ together with a $\tilde c^-  \in \Gamma ^- (  \tilde \sigma_-   )$ be given  with some $q \in \Bbb N$ by
$$
 \Phi(c^-  \sigma_-  d^- \alpha_- ^{2L + 1}  )   = \tilde c^- \tilde \sigma^- \tilde d^- \tilde \alpha_- ^q. 
$$
Let $\tilde b^-$ be the suffix of length $Q+L$ of $\tilde c^-$. Then (1) will hold.

$\sigma _-^+$ and $d_+^-$ have the analogous meaning, and $I_+, \sigma _+, d^+, \sigma _+^-,d^+_-$ have the symmetric meaning.

We prove that
 $$
 \tilde\Xi_-(\tilde X) \ne \emptyset, \tag 2
 $$ 
and that
$$
 \{   \tilde \sigma_-, \tilde d^-:    \tilde\sigma_- \in \Sigma_-( \tilde X) ,d^- \in \Cal D (  \tilde \sigma_-  ,   \tilde\alpha_-   )  \} \setminus  \Xi_-( \tilde X)
\neq \emptyset. \tag 3
$$

Let
$$
\tilde \sigma _- \in \Sigma_- (\tilde X), \quad \tilde d^- \in \Cal D ( \tilde \sigma _-  , \tilde  \alpha_-), \quad \tilde b^-\in \Gamma^-_{Q+L}( \tilde \sigma _-  ).
$$
To prove (2) and (3) we show that 
$$
\sigma_-(\tilde b_- \tilde \sigma _-   )d^- ( \tilde b^-\tilde \sigma _- \tilde d^-     ) \in \Xi _-(X), \tag 4
$$
if and only if
$$
\tilde \sigma _- \tilde d^-  \in \Xi _-(\tilde X). \tag 5
$$
Assume (4), let
$$
\tilde \sigma^-_+ \in \Sigma^-_+ (\tilde X), \quad \tilde d_-^+ \in \Cal D(\tilde\alpha_-,  \tilde \sigma_+^- ),
$$
and let
$$
k_- > \max \  ( 2L, K_-(X)). \tag 6
$$
(4) and (6) imply that for $  \tilde b^+_-   \in \Gamma^+_{L+Q}(   \tilde\sigma^-_+ )$ the word
$$
\multline
\sigma_-(\tilde b_- \tilde \sigma _-   )d^- ( \tilde b^-\tilde \sigma _- \tilde d^-     ) \\
\alpha_-^
{  \ell ( \tilde d^-  ) -\ell  ( d^- ( \tilde b^-\tilde \sigma _- \tilde d^-     )) + I_-( \tilde b_- \tilde \sigma _-    ) + k_-   + I_+( \tilde\sigma^-_+ \tilde b^+_-  )          -  \ell ( d_-^+(   \tilde d_-^+  \tilde\sigma^-_+ \tilde b^+_-  ))              +       \ell  (  \tilde d^+ _-  )   }\\
 d_-^+(   \tilde d_-^+  \tilde\sigma^-_+ \tilde b^+_-  )\sigma^-_+(   \tilde\sigma^-_+ \tilde b^+_-  ) 
 \endmultline
$$
is in
$ \Cal L (X$).
It follows that
$$
 \tilde b^-\tilde \sigma _- \tilde d^-   \tilde \alpha_-^{k_- }  \tilde d_-^+  \tilde\sigma^-_+ \in \Cal L (\tilde X),
$$
and (5) is shown.

Assume (5), and let
$$
 \sigma_+^- \in \Sigma_-^+ (X), \quad d_-^+ \in \Cal D(\alpha_-,  \sigma_+^- ).
$$
Choose
$$
\tilde \sigma_+ ^-\in \Sigma_-^+(\tilde X),\quad  \tilde d^+_- \in \Cal D(  \tilde \alpha_- ,\tilde \sigma_+ ^-  ), \quad \tilde b^+_ -\in \Gamma^+_{Q+ L}(   \tilde \sigma_+ ^-),
$$
such that
$$
\sigma_+ ^- =  \sigma_+ ^- ( \tilde \sigma_+ ^-\tilde b^+_ -  ) , \quad d _-^+  = d_-^+ (  \tilde d^+_-  \tilde \sigma_+ ^-\tilde b^+_ -  ),
$$
and let
$$
\tilde k_- > \max \ (2L, K_-( \tilde X)). \tag 7
$$
(5) and (7) imply that
$$
\tilde\sigma_-\tilde d ^- \tilde\alpha_- 
^{     \ell (   d^-( \tilde b^- \tilde\sigma_- \tilde d^-   )       - \ell (\tilde d ^-) - I_- ( \tilde b^- \tilde\sigma_-  ) + \tilde k_- - I_+( \tilde \sigma_+ ^- \tilde b^+_ -   )         + \ell( d^+_-  ) - \ell (\tilde  d^+_-  )             }
 \tilde d^+_- \tilde \sigma_+ ^- \in \Cal L (\tilde X).
$$
It follows that
$$
\sigma_-(\tilde b^- \tilde\sigma_- )d^-( \tilde b^- \tilde\sigma_- \tilde d^-   ) \alpha_-^
{\tilde k_- } d^+_-\sigma_+^- \in  \Cal L ( X),
$$
and (4) is shown. The proof that
 $$
 \tilde\Xi_+(\tilde X) \ne \emptyset,
\quad  \{  \tilde d^+  \tilde \sigma_+:    \tilde\sigma_+ \in \Sigma_-( \tilde X) ,d^+ \in \Cal D ( \tilde\alpha_+   , \tilde \sigma_+     )  \} \setminus  \Xi_+( \tilde X)
\neq \emptyset,
$$
is symmetric.

We set
$$
\mu_- = \max_{\sigma_- \in \Sigma_-(X). d^- \in \Cal D(\sigma_-  , \alpha_-  )} \ell (d^-),
$$
$\mu_+$ has  the symmetric meaning.

We prove that  $\tilde X$ satisfies condition $(R^-)$. For this let
$$
\tilde \sigma_-^+ \in \Sigma  (\tilde X), d_+^- \in \Cal D ( \tilde \sigma_-^+ , \alpha_+),
$$
such that
$$
\ell( d_+^-   ) > \max ( R^- (X) + L, \mu^- + 2L), \tag 8
$$
and let $\tilde b_+^- \in \Gamma^-_{L + Q}( \tilde \sigma_-^+  )$. By (8) 
$$
\ell (d_+^-( \tilde b_+^- \tilde \sigma_-^+ d_+^- ))  > R^-(X),
$$
and it follows that 
$$
\sigma_-^+(  \tilde b_+^- \tilde \sigma_-^+  ) \in  \Sigma_-(X),
$$
and that there is a 
$$
d^- \in \Cal D (  \sigma_-^+(  \tilde b_+^- \tilde \sigma_-^+  ) , \alpha_-)
$$
such that
$$
d^-_+ (  \tilde b_+^- 
\tilde \sigma_-^+\tilde d^-_+  ) =
d^-  \alpha_-^{\ell (d^-_+ (  \tilde b_+^- 
\tilde \sigma_-^+\tilde d^-_+  ) ) - \ell ( d^-  ) - \ell (c_ X)}c_ X.
$$
It follows that
$$
\tilde \sigma_-^+ \in \Sigma_-(\tilde X),
$$
and with $\tilde d^-\in \Cal D( \tilde \sigma_-^+  ,\tilde \alpha_-)$ such that 
$$
d^- = d^-( \tilde b_+^-\tilde \sigma_-^+  \tilde d^- ),
$$
one has
$$
\tilde d^-_+ =  \tilde d^- \tilde \alpha^{\ell (\tilde d^-_+ ) - \ell (\tilde d^-  ) - \ell (c_{\tilde X})}  .
$$
The proof 
that  $\tilde X$ satisfies condition $(R^+)$ is symmetric.

We prove that $\tilde X$  satisfies conditions $(R^-_\Xi)$.  
For this let $\tilde \sigma_-^+ \tilde d_-^+ \in \Xi _-^+  ( X),$
$$
\ell (  \tilde d_-^+  ) > \max (\mu^-(X) +\ell (c_X) + 2L, R^-_\Xi (X) + L), \tag 9
$$
and let $\tilde b_-^+ \in \Gamma^-_{L + Q}( \tilde \sigma_-^+  )$. Then one has that
$
\sigma_-^+(   \tilde b_-^+ \tilde \sigma_-^+ ) d^-_+(  \tilde b_-^+   \tilde \sigma_-^+ \tilde d_-^+  ) \in \Xi _- ^+ (X).
$
By (9),
$
\ell (  d^-_+(  \tilde b_-^+   \tilde \sigma_-^+ \tilde d_+^-  ) ) > R_\Xi^- (X), 
$
and it follows that $ \sigma_-^+(   \tilde b_-^+ \tilde \sigma_-^+ )\in \Sigma_-(X),$ 
and that there is a $d_-\in \Cal D(  \sigma_-^+(   \tilde b_-^+ \tilde \sigma_-^+ ), \alpha_-)$ such that 
$  
\sigma_-^+(   \tilde b_-^+ \tilde \sigma_-^+ )d_-\in \Xi_-(X),
$
and $$
 d^-_+(  \tilde b_-^+   \tilde \sigma_-^+ \tilde d_-^+  )=  d_-\alpha_-^{\ell ( d^-_+(  \tilde b_-^+   \tilde \sigma_-^+ \tilde d_-^+  ) ) - \ell ( d_- ) -\ell (c_X)} c_X.
$$
By (9)
$$
\ell ( d^-_+(  \tilde b_-^+   \tilde \sigma_-^+ \tilde d_-^+  ) ) - \ell ( d_- ) -\ell (c_X)> 2L. \tag 10
$$
It follows that $\tilde \sigma_-^+ \in \Sigma_-(\tilde X)$, and that one has  with $\tilde d_- \in \Cal D(\tilde \sigma_-^+    , \tilde \alpha_-)$ such that 
$$
d_ -= d_-(  \tilde b_-^+   \tilde \sigma_-^+ \tilde d_-    )
$$
that
$$
 \tilde \sigma_-^+ \tilde d_-   \in \Xi_-(\tilde X),
$$
and as a consequence of  (10) there
$$
 \tilde d_-^+  = \tilde d_- \tilde \alpha_- ^{\ell(  \tilde d_-^+  ) - \ell (  \tilde d^-) - \ell (c_{\tilde X})}c_{\tilde X}.
$$
The proof that $\tilde X$ satisfies condition $(R^+_\Xi)$ is symmetric.

Let $I, J_-, J_+, \Delta _-,   \Delta _+$ be parameters for $X$, and let

$$
\tilde I \geq \max (I,
\mu_- + \mu_+  +\ell (c_X) + 4L +2,  ,
R^-_\Xi  +  R^+_\Xi   + 2Q + 2L   ).  \tag 11
$$
For $ \tilde I _- , \tilde I _+\in \Bbb Z$, such that
$$
\tilde I _+ -  \tilde I _- > \tilde I 
$$
and for
$$
\tilde b \in  \tilde X_{[\tilde I_-, \tilde I_+]}
$$
every choice of $ \tilde b^- \in \Gamma_{Q+L}( \tilde b),  \tilde b^+\in \Gamma_{Q+L}( \tilde b)$ determines $I_-, I_+ \in \Bbb Z$, such that 
$$
\tilde I _- - L - Q \leq  I _- \leq \tilde I _-, \quad  \tilde I _+ \leq   I _+ \leq \tilde I _- + L + Q ,
$$
and such that the word
$$
b=(\Phi (  \tilde b^-\tilde b\tilde b^+ ))_{[I_-, I_+]}
$$
is in $\Cal B(X)$.
Under the hypothesis (11), if here
$$
b\in   \Cal B^{(X)} (\Xi, \alpha_ - , \alpha_ +)\quad
( \ b \in \Cal B^{(X)} ( \alpha_ - , \alpha_ +, \Xi) \ ),
$$
then, applying conditions $(R^-_\Xi)$ and $(R^+_\Xi)$, one confirms that
$$
\tilde b \in 
\Cal B^{(\tilde X)} (\Xi, \tilde\alpha_ - )\cup
\Cal B^{(\tilde X)} (\Xi, \tilde\alpha_ - , \tilde\alpha_ +)
\cup \Cal B^{(\tilde X)} (\Xi, \tilde\alpha_ +)
$$
$$
( \ \tilde b\in\Cal B^{(\tilde X)} (\tilde \alpha_ - ,\tilde \alpha_ +, \Xi)  \cup \Cal B^{(\tilde X)} (\tilde \alpha_ +, \Xi)\cup \Cal B^{(\tilde X)} (\tilde \alpha_ - , \Xi)  \ ),
$$ 
if here
$$
b\in \Cal B^{(X)} (\Xi, \alpha_ - ) \quad
( \ b\in \Cal B^{(X)} (\alpha_ +, \Xi) \ ),
$$
then, applying conditions $(  R^-_\Xi)$ and $(R^+) $ ($(R^-)$ and $(   R^+_\Xi)$), one confirms that
$$
\tilde b\in 
\Cal B^{(\tilde X)} (\Xi, \tilde\alpha_ - )  \cup \Cal B^{(\tilde X)} (\Xi, \tilde\alpha_ - , \tilde\alpha_ +) 
$$
$$
( \ \tilde b \in  \Cal B^{(\tilde X)} \tilde \alpha_ +, \Xi)  \cup\Cal B^{(\tilde X)} (\tilde \alpha_ - ,\tilde \alpha_ +, \Xi) \ ),
$$ 
and if here
$$
b\in \Cal B^{(X)} (\Xi, \alpha_+  ) \quad
( \ b\in \Cal B^{(X)} (\alpha_ -, \Xi) \ ),
$$
then, applying condition $(  R^-_\Xi)$ and $(R^+) $ ($(R^-)$ and $(   R^+_\Xi)$), one confirms that
$$
\tilde b \in 
\Cal B^{(\tilde X)} (\Xi, \tilde\alpha_ +)  \cup \Cal B^{(\tilde X)} (\Xi, \tilde\alpha_ - , \tilde\alpha_ +) 
$$
$$
( \ \tilde b\in  \Cal B^{(\tilde X)} \tilde \alpha_ -, \Xi)  \cup\Cal B^{(\tilde X)} (\tilde \alpha_ - ,\tilde \alpha_ +, \Xi) \ ).
$$ 
We define 
$ H_-, H_+ \in \Bbb Z_+$ 
by
$$
\tilde\Phi (c_{\tilde X} )= \alpha_-^{H_-}c_X \alpha_+^{H_+},
$$
and we set
$$
\tilde J_- = J_-  + H_- , \quad  \tilde J_+ = J_+  + H_+ .
$$
Under the hypothesis (11), if here
$$
b \in  \Cal B_{counter}^{(X)}(\Delta^-, J_-, J_+,\Delta^+) 
$$
then, setting
$$
\align
\tilde{\Delta}_-(\tilde{\sigma}_-\tilde{d}^-)
& = \bigcup_{\tilde{b}^- \in \Gamma^-_{Q +L}(\tilde{\sigma}_-)}
 \Delta_-(\sigma_-(\tilde{b}^-\tilde{\sigma}_-)
    d^-(\tilde{b}^-\tilde{\sigma}_-\tilde{d}^-)) + \ell({\tilde{d}^-}) \\
& \qquad 
 -\ell(d^-(\tilde{b}^-\tilde{\sigma}_-\tilde{d}^-)) 
 + I_-(\tilde{b}^-\tilde{\sigma}_-)),\\
& \qquad \qquad
 \tilde{\sigma}_- \in \tilde{\Sigma}_-(\tilde{X}),
 \tilde{d}^- \in {\Cal D}(\tilde{\sigma}_-, \tilde{\alpha}_-),  \tilde{\sigma}_- \tilde{d}^- \notin \Xi_-(\tilde{X}),
 \endalign
 $$
 $$
\align
\tilde{\Delta}_+(\tilde{d}^+\tilde{\sigma}_+)
& = \bigcup_{\tilde{b}^+ \in \Gamma^+_{L+Q}(\tilde{\sigma}_+)}
 I_+(\tilde{\sigma}_+\tilde{b}^+)-\ell(d^+(\tilde{d}^+\tilde{\sigma}_+\tilde{b}^+)) + \ell(\tilde{d}^+) \\
& \qquad + \Delta_+(d^+(\tilde{d}^+\tilde{\sigma}_+\tilde{b}^+)\sigma_+(\tilde{\sigma}_+\tilde{b}^+)),\\
& \qquad \qquad
 \tilde{\sigma}_+ \tilde{d}^+\notin \Xi_+(\tilde{X}).
\endalign
$$
one confirms that 
$$
 \tilde b\in  \Cal B_{counter}^{(\tilde X)}(\tilde\Delta^-, \tilde J_-, \tilde J_+,\tilde\Delta^+) .
$$
It can also be proved that
$$
 \Cal B_{counter}^{(\tilde X)}(\tilde\Delta^-, \tilde J_-, \tilde J_+,\tilde\Delta^+) 
\subset \Cal L( \tilde X).
$$
This means that $\tilde I, \tilde\Delta^-, \tilde J_-, \tilde J_+,\tilde\Delta^+$ are parameters of the Lind-Marcus shift $\tilde X$.
\qed
\enddemo

\proclaim{Theorem 4.2} 
A subshift  that is topologically conjugate to a Lind-Marcus one-counter shift has an $n$-block system that is a Lind-Marcus one-counter shift.
\endproclaim
\demo{Proof}
Let $X_{\circ} \subset \Sigma_{\circ} ^{\Bbb Z}  $ be a subshift that is topologically conjugate to Lind-Markov one-counter shift. By Lemma 2.3 of \cite {KM} $X_{\circ} $   has an $n$-block system $\tilde X$ such that there is a topological conjugacy $\tilde\varphi : \tilde X \to X $ that is given by a one-block map
$\tilde\varPhi:\tilde\Sigma \to {\Sigma}$,
such that 
$$
\tilde\varPhi^{-1}(\Sigma_{synchro}(X) \subset \Sigma_{synchro}(\tilde X).
$$
Apply Lemma 4.1.
\qed
\enddemo

\bigskip
{\it e-mail}: 

krieger{\@}math.uni-heidelberg.de

\bye